\documentclass[12pt]{amsart}
\usepackage{amsmath,amsfonts,amssymb,amsthm}
\usepackage{rotating}
\def\NN{\mathbb{ N}}

\def\Z{{\mathbb Z}}

\def\inv{\textrm{inv}}
\def\AP{\mathcal{AP}}
\def\SAP{\mathcal{SAP}}
\newcommand\pv{\noindent{\bf Proof : \/}}
\def\qed{\quad\raise -2pt\hbox{\vrule\vbox
to 10pt{\hrule width 4pt \vfill\hrule}\vrule}\medskip}
\newtheorem{thm}{Theorem}
\newtheorem{lem}{Lemma}

\newtheorem{prop}{Proposition}

\newtheorem{cor}{Corollary}
\newcommand\ch{\mathop{\rm ch}}
\def\C{{\mathbb C}}

\title{A $Q$-analog of the Seidel generation of Genocchi numbers}
\begin{document}
\maketitle \centerline{Jiang Zeng$^1$ and Jin Zhou$^2$}
\begin{center} \small $^1$ Institut Girard Desargues,
Universit\'e Claude Bernard (Lyon I)\\
69622 Villeurbanne Cedex, France \\
{\tt zeng@igd.univ-lyon1.fr}\\
\vspace{10pt}
 $^2$ Center for Combinatorics, LPMC,
Nankai University\\
 Tianjin 300071, People's Republic of China\\
{\tt jinjinzhou@hotmail.com} \\
\end{center}

\begin{abstract}
A  new $q$-analog of Genocchi numbers is introduced through a
$q$-analog of Seidel's triangle associated to Genocchi numbers. It
is then shown that these $q$-Genocchi numbers have interesting
combinatorial interpretations in the classical models for Genocchi
numbers such as alternating pistols, alternating permutations, non
intersecting lattice paths and skew Young tableaux.
\end{abstract}

\section{Introduction}
The \emph{Genocchi numbers} $G_{2n}$ can be defined  through their
relation with Bernoulli numbers $G_{2n}=2(2^{2n}-1)B_n$ or by
their exponential generating function~\cite[p. 74-75]{St2}:
$$
\frac{2t}{e^t+1}=t-\frac{t^2}{2!}+\frac{t^4}{ 4!}-3\,
\frac{t^6}{6!} +\cdots +(-1)^nG_{2n}\frac{t^{2n}}{(2n)!}+\cdots.
$$
 However it is
not straightforward from the above definition
 that  $G_{2n}$ should be  \textit{integers}. It was Seidel~\cite{Se} who  first gave  a
  Pascal type triangle for  Genocchi numbers in the nineteenth century.
Recall that the \emph{Seidel triangle} for Genocchi
numbers~\cite{DV,DZ1,Vi} is an array of integers
$(g_{i,j})_{i,j\geq 1}$ such that $g_{1,1}=g_{2,1}=1$ and
\begin{equation}\left\lbrace
\begin{array}{lll}
      g_{2i+1,j}&=&g_{2i+1,j-1}+g_{2i,j},
\quad \hbox{for}\quad j=1,2,\ldots,i+1,\\
g_{2i,j}&=&g_{2i,j+1}+g_{2i-1,j}, \quad \hbox{for}\quad j=i,
i-1,\ldots, 1,
\end{array}
\right.
\end{equation}
where $g_{i,j}=0$ if $j<0$ or $j> \lceil{i/2}\rceil$ by
convention. The first values of $g_{i,j}$ for $1\leq i,j\leq 10$
can be displayed in \emph{Seidel's tiangle for Genocchi numbers}
as follows:
\begin{table}[h]
\begin{tabular}{ccccccccccc|}
\cline{9-11} &&&&&&&&\multicolumn{1}{|c}{155}&\multicolumn{1}{c}{155}&\multicolumn{1}{|c|}{5}\\
\cline{7-8} &&&&&&\multicolumn{1}{|c}{17}&\multicolumn{1}{c}{17}&155&310&\multicolumn{1}{|c|}{4}\\
\cline{5-6}&&&&\multicolumn{1}{|c}{3}&
\multicolumn{1}{c}{3}&\multicolumn{1}{c}{17}&
\multicolumn{1}{c}{34}&138&448&\multicolumn{1}{|c|}{3}\\
\cline{3-4} &&\multicolumn{1}{|c}{1}&\multicolumn{1}{c}{1}&
\multicolumn{1}{c}{3}&\multicolumn{1}{c}{
6}&\multicolumn{1}{c}{14}&\multicolumn{1}{c}{48} &104&552&\multicolumn{1}{|c|}{2}\\
\cline{1-2}\multicolumn{1}{|c}{1}&
\multicolumn{1}{c}{1}&\multicolumn{1}{c}{1}&\multicolumn{1}{c}{2}
&\multicolumn{1}{c}
{2}&\multicolumn{1}{c}{8}&\multicolumn{1}{c}{8}
&\multicolumn{1}{c}{56}&56&608&\multicolumn{1}{|c|}{1}\\
\cline{1-11}  \multicolumn{1}{|c}{1}&\multicolumn{1}{c}{2}&
\multicolumn{1}{c}{3}&\multicolumn{1}{c}{4}&\multicolumn{1}{c}{5}&
\multicolumn{1}{c}{6}&\multicolumn{1}{c}{7}&\multicolumn{1}{c}{8}
&\multicolumn{1}{c}{9}&\multicolumn{1}{c}{10}&\multicolumn{1}{|c|}{$i\setminus j$}\\
\hline
\end{tabular}
\end{table}

The Genocchi numbers $G_{2n}$ and the so-called \emph{median
Genocchi numbers} $H_{2n-1}$ are given  by the
following relations~\cite{DV}:
$$
G_{2n}=g_{2n-1,n},\qquad H_{2n-1}=g_{2n-1,1}.
$$

The purpose of this paper is to show that there is a $q$-analog of
Seidel's algorithm  and the resulted $q$-Genocchi numbers inherit
most of the nice results proved by Dumont-Viennot, Gessel-Viennot
and Dumont-Zeng for ordinary Genocchi numbers~\cite{DV,GV2,DZ2}.

Note that some different $q$-analogs of Genocchi numbers have been
investigated from  both combinatorial and algebraic points of
view~\cite{HZ,Ra}. In particular, Han and Zeng~\cite{HZ} have
found an interesting $q$-analog of Gandhi's algorithm ~\cite{Ga}
by using the $q$-difference operator instead of the difference
operator and proved that the ordinary generating function of these
$q$-Genocchi numbers has a remarkable continued fraction
expansion.

A $q$-Seidel triangle is  an array  $(g_{i,j}(q))_{i,j\geq 1}$ of polynomials in $q$
such that $g_{1,1}(q)=g_{2,1}(q)=1$ and
\begin{equation}\label{eq3}\left\lbrace
\begin{array}{lll}
      g_{2i+1,j}(q)&=&g_{2i+1,j-1}(q)+q^{j-1}g_{2i,j}(q),
\quad \hbox{for}\quad j=1,2,\ldots,i+1,\\
g_{2i,j}(q)&=&g_{2i,j+1}(q)+q^{j-1} g_{2i-1,j}(q), \quad
\hbox{for}\quad j=i, i-1,\ldots, 1,
\end{array}
\right.
\end{equation}
where $g_{i,j}(q)=0$ if $j<0$ or $j> \lceil{i/2}\rceil$
 by
convention. The first values of $g_{i,j}(q)$  are given in Table~\ref{firgij}.
\vspace{-1.5cm}
%%%%%%%%%%%%%%%%%%%%%%%%%%%%%%%%%%%%%%%%%%%%%%%%%%%%%%%%%%%%%%%%%%%%%%
\begin{table}[h]{\scriptsize
\begin{tabular}{cccccccccc|}
\multicolumn{10}{c}{ }\\ \multicolumn{10}{c}{ }\\
\multicolumn{10}{c}{ }\\ \multicolumn{10}{c}{ }\\ \cline{9-10}
&&&&&&&&\multicolumn{1}{|c|}{$1+2q+3q^2+4q^3+4q^4+2q^5+q^6$} &4\\
\cline{7-8} &&&&&&\multicolumn{1}{|c}{$1+q+q^2$}&
\multicolumn{1}{c|}{$q^2+q^3+q^4$}&\multicolumn{1}{|c|}{$1+2q+3q^2+4q^3+4q^4+2q^5+q^6$}
&3\\ \cline{5-6} &&&
&\multicolumn{1}{|c}{$1$}&\multicolumn{1}{c|}{$q$}&
\multicolumn{1}{|c}{$1+q+q^2$}&\multicolumn{1}{c|}{
$q+2q^2+2q^3+q^4$}&\multicolumn{1}{|c|}{$1+2q+3q^2+4q^3+3q^4+q^5$}
&2\\ \cline{3-4} &&
\multicolumn{1}{|c}{1}&\multicolumn{1}{c|}{1}&\multicolumn{1}{|c}{1}&\multicolumn{1}{c|}
{$1+q$}&\multicolumn{1}{|c}{$1+q$}&\multicolumn{1}{c|}{$1+2q+2q^2+2q^3+q^4$}&\multicolumn{1}{|c|}
{$1+2q+2q^2+2q^3+q^4$}&1\\ \cline{3-10}
&&\multicolumn{1}{|c|}{1}&\multicolumn{1}{|c|}{2}&\multicolumn{1}{|c|}{3}&
\multicolumn{1}{|c|}{4}&\multicolumn{1}{|c|}{5}
&\multicolumn{1}{|c|}{6}&\multicolumn{1}{|c|}{7}&$i\setminus j$\\
\cline{3-10} \multicolumn{10}{c}{ }\\
\end{tabular}}
\caption{\label{firgij}$q$-analog of Seidel's triangle
$(g_{i,j}(q))_{i,j\geq 1}$}
\end{table}
%%%%%%%%%%%%%%%%%%%%%%%%%%%%%%%%%%%%%%%%%%%%%%%%%%%%

Define the $q$-Genocchi numbers $G_{2n}(q)$ and
 \emph{$q$-median Genocchi numbers} $H_{2n-1}(q)$ by $G_2(q)=H_1(q)=1$ and
 for all $n\geq 2$~:
\begin{equation}
G_{2n}(q)=g_{2n-1,n}(q),\quad H_{2n-1}(q)=q^{n-2}g_{2n-1,1}(q).
\end{equation}
Thus, the sequences for $G_{2n}(q)$ and $H_{2n-1}(q)$ start with
$1, 1, 1+q+q^2$ and $1,1,q+q^2$, respectively.

This paper is organised as follows.
In sections 2 and 3 we generalize the combinatorial results of Dumont and Viennot~\cite{DV}
by first interpreting  $g_{i,j}(q)$ (and in
particular the two kinds of $q$-Genocchi numbers) in the model of
alternating pistols and then derive the interpret $G_{2n}(q)$ as generating polynomials of
 \emph{alternating permutations}.
 In section~4 we give the $q$-version of the results of
Gessel-Viennot~\cite{GV2} and Dumont-Zeng~\cite{DZ1}. In
section~4, by extending the matrix of $q$-binomial coefficients to
\emph{negative indices} we obtain a $q$-analog of results of
Dumont and Zeng~\cite{DZ2}. Finally, in section 6, we show that
there is a remarkable triangle of $q$-integers containing the two
kinds of $q$-Genocchi numbers and conjecture that the terms of
this triangle refine the classical $q$-secant numbers,
generalizing a result of Dumont-Zeng~\cite{DZ1}.

%%%%%%%%%%%%%%%%%%%%%%%%%%%%%%%%%%%%%%%%%%%%%%%%%%%%%%%%%%%%%%%%%%%%%%%
\section{Alternating pistols}
%%%%%%%%%%%%%%%%%%%%%%%%%%%%%%%%%%%%%%%%%%%%%%%%%%%%%%%%%%%%%%%%%%%%%%%%
An \emph{alternating pistol} (resp. \emph{strict-alternating
pistol}) on $[m]=\{1,\cdots,m\}$ is a mapping  $p:[m] \to [m]$
such that for $i=1,2,\ldots, \lceil m/2\rceil$:
\begin{enumerate}
\item $p(2i)\leq i$ and $p(2i-1)\leq i$,
\item $p(2i-1)\geq p(2i)$ and $p(2i) \leq p(2i+1)$ \hspace{.5cm}(resp. $p(2i)<
p(2i+1)$).
\end{enumerate}
\goodbreak
We can illustrate an alternating pistol on $[m]$ by an array
$(T_{i,j})_{1\leq i,j\leq m}$ with a cross at $(i,j)$ if $p(i)=j$.
For example, the alternating pistol  $p=p(1)p(2)\ldots p(8)=11211143$ can be
 illustrated as in Figure~\ref{alterp}.\\

%%%%%%%%%%%%%%%%%%%%%%%%%%%%%%%%%%%%%%%%%%%%%%%%%%%%%%%%%%%%%%%%%%%%%%%%
\begin{figure}[h]{\scriptsize
\begin{tabular}{ccccccccc|}
\cline{7-9} &&&&&&\multicolumn{1}{|c|}{X}
&\multicolumn{1}{|c|}{}&4\\
\cline{5-8}&&&&\multicolumn{1}{|c|}{}&
\multicolumn{1}{|c|}{}&\multicolumn{1}{|c|}{}&
\multicolumn{1}{|c|}{X}&3\\
\cline{3-8}&&\multicolumn{1}{|c|}{X}&\multicolumn{1}{|c|}{}&
\multicolumn{1}{|c|}{}&\multicolumn{1}{|c|}{
}&\multicolumn{1}{|c|}{}&\multicolumn{1}{|c|}{}&2\\
\cline{1-8}\multicolumn{1}{|c|}{X}&
\multicolumn{1}{|c|}{X}&\multicolumn{1}{|c|}{}&\multicolumn{1}{|c|}{X}
&\multicolumn{1}{|c|}
{X}&\multicolumn{1}{|c|}{X}&\multicolumn{1}{|c|}{}&\multicolumn{1}{|c|}
{}&1\\
 \cline{1-9} \multicolumn{1}{|c}{1}&\multicolumn{1}{c}{2}&
\multicolumn{1}{c}{3}&\multicolumn{1}{c}{4}&\multicolumn{1}{c}{5}&
\multicolumn{1}{c}{6}&\multicolumn{1}{c}{7}&\multicolumn{1}{c}{8}&
\multicolumn{1}{|c|}{$i\setminus j$}\\ \cline{1-9}\\
\end{tabular}
\caption{\label{alterp}An alternating pistol $p=11211143$}}
\end{figure}
%%%%%%%%%%%%%%%%%%%%%%%%%%%%%%%%%%%%%%%%%%%%%%%%%%%%%%%%%%%%%%%%%%%%%%%%%%%
For all $i\geq 1$ and $1\leq j\leq \lceil i/2\rceil$, let
$\AP_{i,j}$ (resp. $\SAP_{i,j}$) be the set of alternating pistols
$p$ (resp. strict-alternating pistols) on $[i]$ such that
$p(i)=j$. Dumont and Viennot~\cite{DV} proved that the entry
$g_{i,j}$ of  Seidel's triangle is the cardinality of $\AP_{i,j}$.
Hence $G_{2n}$ (resp. $H_{2n+1}$) is the number of alternating
pistols (resp. strict alternating pistols) on $[2n]$.

To obtain a $q$-version of Dumont-Viennot's result, we define the
\emph{charge} of a pistol $p$ by
$$
{\ch} (p)= (p_1-1)+(p_2-1)+\cdots +(p_{m}-1).
$$
In other words the charge of a pistol $p$ amounts to the number of
cells below its crosses. For example, the charge of the pistol
 in Figure~\ref{alterp} is $ch(p)=1+3+2=6$.

\begin{prop}
For $i\geq 1$ and $1\leq j\leq \lceil i/2\rceil$,
$g_{i,j}(q)$ is the generating function of alternating pistols $p$
on $[i]$ such that $p(i)=j$, with respect to
the charge, i.e.,
$$
g_{i,j}(q)=\sum_{p\in \mathcal{AP}_{i,j}}q^{{\ch}(p)-j+1}.
$$
\end{prop}
\begin{pv}
We proceed by double inductions on $i$ and $j$, where $1\leq j\leq
\lceil i/2\rceil$:
\begin{itemize}
\item If $i=1$, then $p(1)=1$ and
${\ch}(p)=0$, so $g_{1,1}(q)=1$,
\item Let $p\in\mathcal{AP}_{2k+1,j}$ and suppose the recurrence is
true for all elements of $\mathcal{AP}_{2k'+1,j'}$ with $k'<k$, or
$k'=k$ and $j'<j$.
\begin{enumerate}
\item If $j>p(2k)$, let $p'\in \mathcal{AP}_{2k+1,j-1}$ such that
$p$ and $p'$ have the same restrictions to $[2k]$. Then
${\ch}(p)={\ch}(p')$,
\item If $j=p(2k)$ then the charge of the restriction of $p$ to
$[2k]$ is ${\ch}(p)-j+1$.
\end{enumerate}
Summing over all elements of $\mathcal{AP}_{2k+1,j}$, we obtain
the first equation of (\ref{eq3}).
\item  Let $p\in\mathcal{AP}_{2k,j}$
and suppose the recurrence true for all elements of
$\mathcal{AP}_{2k',j'}$ with $k'<k$, or $k'=k$ and $j'>j$.
\begin{enumerate}
\item If $j<p(2k-1)$, let $p'\in \mathcal{AP}_{2k,j+1}$
such that $p$ and $p'$ have same restrictions to $[2k-1]$. Then
${\ch}(p)={\ch}(p')$.
\item If $j=p(2k-1)$ then the charge of the restriction of $p$ to
$[2k-1]$ is ${\ch}(p)-j+1$.
\end{enumerate}
\end{itemize}
Summing over all elements of $\mathcal{AP}_{2k,j}$, we obtain the
second equation of (\ref{eq3}). \qed
\end{pv}
%%%%%%%%%%%%%%%%%%%%%%%%%%%%%%%%%%%%%%%%%%%%%%%%%%%%%%%%%%%%%%%%%%%%%%%%%%%%%%%%

\begin{table}{\scriptsize
\begin{tabular}{ccccccccc|} \multicolumn{9}{c}{
}\\\multicolumn{9}{c}{ }\\ \cline{7-9}
&&&&&&\multicolumn{1}{|c}{$q^2+2q^3+2q^4+2q^5+q^6$}&\multicolumn{1}
{c|}{$q^5+2q^6+2q^7+2q^8+q^9$}&4\\ \cline{5-6}
&&&&\multicolumn{1}{|c}{$q+q^2$}&\multicolumn{1}{c|}{
$q^3+q^4$}&\multicolumn{1}{|c}{$q^2+2q^3+2q^4+q^5$}&\multicolumn{1}{c|}
{$q^4+3q^5+4q^6+3q^7+2q^8+q^9$}&3\\ \cline{3-4} &
&\multicolumn{1}{|c}{$1$}&\multicolumn{1}{c|}{$q$}&
\multicolumn{1}{|c}{$q$}&\multicolumn{1}{c|}{
$q^2+q^3+q^4$}&\multicolumn{1}{|c}{$q^2+q^3+q^4$}&\multicolumn{1}{c|}
{$q^3+2q^4+4q^5+4q^6+3q^7+2q^8+q^9$}&2\\ \cline{1-2}
\multicolumn{1}{|c}{1}&\multicolumn{1}{c|}{1}&\multicolumn{1}{|c}{0}&\multicolumn{1}{c|}
{$q$}&\multicolumn{1}{|c}{0}&\multicolumn{1}{c|}{$q^2+q^3+q^4$}&\multicolumn{1}{|c}{0}&
\multicolumn{1}{c|}{$q^3+2q^4+4q^5+4q^6+3q^7+2q^8+q^9$}&1\\
\cline{1-9}
\multicolumn{1}{|c|}{1}&\multicolumn{1}{|c|}{2}&\multicolumn{1}{|c|}{3}
&\multicolumn{1}{|c|}{4}&\multicolumn{1}{|c|}{5}
&\multicolumn{1}{|c|}{6}&\multicolumn{1}{|c|}{7}&\multicolumn{1}{|c|}{8}&$i\setminus
j$\\ \cline{1-9} \multicolumn{9}{c}{ }
\\
\end{tabular}}
\caption{First values of $h_{i,j}(q)$}\label{firhij}
\end{table}
%%%%%%%%%%%%%%%%%%%%%%%%%%%%%%%%%%%%%%%%%%%%%%%%%%%%%%%%%%%%%%%%%%%%%%
In order to interpret the $q$-median Genocchi numbers
$H_{2n-1}(q)$, it is convenient to introduce another array
$(h_{i,j}(q))_{i,j\geq 1}$ of polynomials in $q$ such that
$h_{1,1}(q)=h_{2,1}(q)=1$, $h_{2i+1,1}(q)=0$ and
\begin{equation}\label{eq4}\left\lbrace
\begin{array}{lll}
         h_{2i+1,j}(q)&=&h_{2i+1,j-1}(q)+q^{j-2}h_{2i,j-1}(q),\\
h_{2i,j}(q)&=&h_{2i,j+1}(q)+q^{j-1}h_{2i-1,j}(q),
\end{array}\right.
\end{equation}
where by convention $h_{i,j}(q)=0$ if $j<0$ or
$j>\lceil{i/2}\rceil$. The first values of $h_{i,j}(q)$ are given
in Table~\ref{firhij}. Similarly we can prove the following:
\begin{prop}
For all $i\geq 1$ and $1\leq j\leq \lceil i/2\rceil$, we have
$$
h_{i,j}(q)=\sum_{\sigma\in \SAP_{i,j}}q^{\ch(\sigma)-j+1}.
$$
\end{prop}
Notice that
$$
G_{2n+2}(q)=g_{2n+1,n+1}(q)=\sum_{1\leq k\leq n}q^{k-1}g_{2n,k}(q),
$$
and since $h_{2n-1,n}(q)=q^{n-2}g_{2n-1,1}(q)$, we have also
$$
H_{2n+1}(q)=h_{2n+1,n+1}(q)=\sum_{1\leq k\leq
n}q^{k-1}h_{2n,k}(q).
$$
The above observations and propositions infer immediately  the following result.
\begin{prop}
For all $n\geq 1$, the $q$-Genocchi number $G_{2n+2}(q)$ (resp.
$q$-medians Gennochi numbers $H_{2n+1}(q)$) is the generating
function of alternating pistols (resp. strict alternating pistols)
on $[2n]$ with respect to the statistics charge, i.e.,
$$
G_{2n+2}(q)=\sum_{p\in \AP_{2n}}q^{\ch p},\qquad H_{2n+1}(q)=\sum_{p\in \SAP_{2n}}q^{\ch p}.
$$
\end{prop}

Dumont and Viennot~\cite[Section 3]{DV} also gave a combinatorial
interpretation of Genocchi numbers with alternating permutations.
In the next section we show that one can translate the statistics
\emph{charge} through all the bijections involved in their proof
and interpret the $q$-Genocchi numbers as a $q$-counting of
alternating permutations.
\section{Alternating permutations}
For any  $\sigma\in S_n$ and $i\in [n]$, the \emph{inversion
table} of $\sigma$ is a mapping $f_{\sigma}: [n]\to [0,n-1]$
defined by:
\begin{center}
$\forall i\in [n], \quad f_{\sigma}(i)$ is the number of indices
$j$ such that $j<i$ and $\sigma(j)<\sigma(i)$.
\end{center}
The mapping $f_\sigma$ is  an \emph{subexceedant function} on
$[n]$, that is a mapping $f_\sigma: [n]\to [0,n-1]$ such that
$0\leq f_\sigma(i)<i$ for every $i\in [n]$.  It is
well-known~\cite[p. 21]{St1} that the correpondance $\ell:
\sigma\mapsto I_\sigma$ is a bijection between the set of
permutations of $[n]$ and  the set of subexceedant functions on
$[n]$. Note that in \cite{St1} the
 \emph{inversion table} of $\sigma$ is the mapping
$I_\sigma: [n]\to [n-1]$ defined by $I_\sigma(i)=i-1-f_\sigma(i)$
for all $i\in [n]$ and the inversion number of a permutation of
$\sigma$ is defined as the following:
\begin{equation}\label{eq:inv}
 \inv \sigma=\sum_{i=1}^n(i-1-f_\sigma(i))=\frac{n(n-1)}{2}-\sum_{i=1}^nf_\sigma(i).
\end{equation}

For example, let $\sigma=8\,3\,9\,4\,5\,1\,6\,2\,7\in S_9$, then
the inversion table is  $f_{\sigma}=002120416$ and the inversion
number is  $\inv \sigma=20$.

A permutation $\sigma$ of $[2n+1]$ is said to be
\emph{alternating} if:
$$
\forall i\in [n],\quad  \sigma(2i-1)>\sigma(2i)\quad
\textrm{and}\quad  \sigma(2i)<\sigma(2i+1).
$$

Let ${\mathcal F}_{2n+1}$ be the set of alternating permutations on $[2n+1]$ with even
inversion table.

%%%%%%%%%%%%%%%%%%%%%%%%%%%%%%%%%%%%%%%%%%%%%%%%%%%%%%%%%%%%%
\begin{prop}
The $q$-Genocchi number $G_{2n+2}(q^2)$ is the generating function
of ${\mathcal F}_{2n+1}$ with respect to $\inv -n$, i.e.,
$$G_{2n+2}(q)=\sum_{\sigma\in {\mathcal F}_{2n+1}}
q^{\frac{1}{2}(\inv \;\sigma-n)}.
$$
\end{prop}
\pv As in \cite{DV}, we define the mapping $\alpha: p\mapsto p'$
from $\AP_{2n}$ to $\AP_{2n+1}$ by
$$
 p'(1)=1,\quad p'(2i)=i+1-p(2i-1),
\quad p'(2i+1)=i+2-p(2i),\quad \forall i\in [n].
$$
Note that ${\ch}(p')=n^2-{\ch}(p)$. Then we can construct an even
subexceedant function $\phi(p')=f$ on $[2n+1]$ by the following
$$
f(i)=2(p'(i)-1), \quad \forall i\in [2n+1].
$$
Let $\sigma=\ell^{-1}(f)$ be the permutation whose inversion table
is$f$, it is easily verified (cf. ~\cite{DV}) that $p$ is an
alternating pistol on $[2n]$ if and only if $\sigma$ is an
alternating permutation $[2n+1]$. Finally, it follows from
\eqref{eq:inv} that
 $$
{\ch}(p)=\frac{1}{2}(\textrm{inv}\sigma-n).
$$
For example, for the alternating pistol $p=11211143\in \AP_8$ in
Figure~\ref{alterp}, we have $p'=112133413\in \AP_9$,
$f=002044604$ and $\sigma=436287915\in \mathcal{F}_9$. \qed

%%%%%%%%%%%%%%%%%%%%%%%%%%%%%%%%%%%%%%%%%%%%%%%%%%%%%%%%%%%%%%%%%%%%%%%%%%%%
\section{Non intersecting lattice paths}
%%%%%%%%%%%%%%%%%%%%%%%%%%%%%%%%%%%%%%%%%%%%%%%%%%%%%%%%%%%%%%%%%%%%%%%%%%%%
The $q$-shifted  factorials $(x;q)_n$ are defined by
$$
(x;q)_n=(1-x)(1-xq)\ldots (1-xq^{n-1}),\qquad \forall n\geq 0.
$$
They can be used to define the $q$-binomial coefficients ${m\brack
n}_q$ as
$$
{m\brack n}_q=\frac{(q^{m-n+1};q)_n}{(q;q)_n}\qquad \forall m\in
\Z\quad \textrm{and}\quad n\in\NN.
$$
Let  $G_q^{-1}=((-1)^{i-j}c_{i,\,j}(q))_{i,\,j\geq 1}$ be the
inverse matrix of
\begin{equation}
G_q=\left({i\brack 2i-2j}_q\,q^{(i-j-1)(i-j)}\right)_{i,\,j\geq
1}.
\end{equation}

The first values of $c_{i,j}(q)$ are given in Table~\ref{fircij}.
%%%%%%%%%%%%%%%%%%%%%%%%%%%%%%%%%%%%%%%%%%%%%%%%%%%%%%%%%%%%%%%%%%%%%%%%%%%%%%%%%%%%%%%%%%
\begin{table}[ht]{\scriptsize
 \begin {tabular}{c|c|c|c|c}
$i\setminus j$&1&2&3&4\\ \hline
1&1&0&0&0\\
2&1&1&0&0\\
3&${q}^{2}+q+1$&${q}^{2}+q+1$&1&0\\
4&$q^6+2q^5+4q^4+4q^3+3q^2+2q+1$&$q^6+2q^5+4q^4+4q^3+3q^2+2q+1$&$\left({q}^{2}+q+1\right
)\left ({q}^{2}+1\right)$&1\\ \multicolumn{5}{c}{ }\\
\multicolumn{5}{c}{ }\\
\end{tabular}}
\caption{First values of $c_{i,j}(q)$ }\label{fircij}
\end{table}

%%%%%%%%%%%%%%%%%%%%%%%%%%%%%%%%%%%%%%%%%%%%%%%%%%%%%%%%%%%%%%%%%%%%%%% We first show that
$c_{k,\,l}(q)$ is a polynomial in $q$ with non negative integrer
coefficients using Gessel-Viennot's theory~\cite{GV1,GV2}.

Let $A$ and  $B$ be two points in the plan $\Pi=\NN\times \NN$ of
coordinates $(a,b)$ and $(c,d)$, respectively. A \emph{lattice
path} from $A$ to $B$ is a sequence of points $((x_i,\,
y_i))_{0\leq i\leq k}$ such that $(x_0,\, y_0)=(a,b)$, $(x_k,\,
y_k)=(c,d)$ and each step is either \emph{east}  or \emph{north},
i.e., $x_i-x_{i-1}=1$ and $y_i-y_{i-1}=0$ or $x_i-x_{i-1}=0$ and
$y_i-y_{i-1}=-1$ for $1\leq i\leq k$. Clearly there is a path from
$A$ to $B$ if and only if $a\leq c$ and $b\geq d$.
%%%%%%%%%%%%%%%%%%%%%%%%%%%%%%%%%%%%%%%%%%%%%%%%%%%%%%

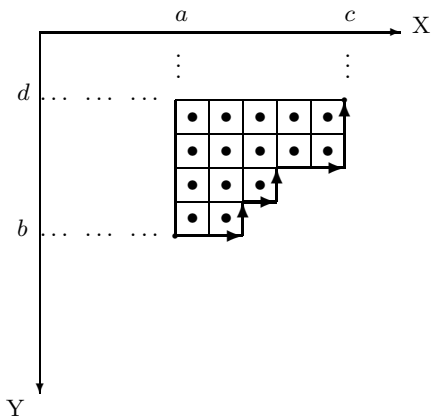
\begin{figure}[ht]{\scriptsize\setlength{\unitlength}{0.15cm}
\begin{picture}(40,40)(-10,-28) \put(1,-2){$\bullet$}\put(4,-2){$\bullet$}\put(7,-2){$\bullet$}
\put(10,-2){$\bullet$}\put(13,-2){$\bullet$}
\put(1,-5){$\bullet$}\put(4,-5){$\bullet$}\put(7,-5){$\bullet$}
\put(1,-8){$\bullet$}\put(4,-8){$\bullet$}\put(7,-8){$\bullet$}
\put(1,-11){$\bullet$}\put(4,-11){$\bullet$}
\put(10,-5){$\bullet$}\put(13,-5){$\bullet$}
\put(0,0){\thinlines\line(1,0){15}}
 \put(3,0){\thinlines\line(0,-1){12}} \put(6, 0){\thinlines\line(0,-1){12}}
\put(9, 0){\thinlines\line(0,-1){9}} \put(12,
0){\thinlines\line(0,-1){6}} \put(15, 0){\thinlines\line(0,-1){6}}
\put(15,0){\circle*{.5}}\put(-14,0){$d$} \put(15,7){$c$}
\put(0,-12){\circle*{.5}}\put(-14,-12){$b$}\put(0,7){$a$}
\put(6,-12){\thicklines\thicklines\vector(0,1){3}}
\put(0,-12){\thicklines\thicklines\vector(1,0){6}}
\put(6,-9){\thicklines\thicklines\vector(1,0){3}}
\put(9,-9){\thicklines\thicklines\vector(0,1){3}}
\put(9,-6){\thicklines\thicklines\vector(1,0){6}}
\put(15,-6){\thicklines\thicklines\vector(0,1){6}}
%\put(-4,-12){$\_$}\put(-1,-12){$\_$}\put(-2,-12){$\_$}\put(-3,-12){$\_$}
%\put(-5,-12){$\_$}\put(-6,-12){$\_$}\put(-7,-12){$\_$}\put(-8,-12){$\_$}
%\put(-10,-12){$\_$}\put(-9,-12){$\_$}\put(-11,-12){$\_$}\put(-12,-12){$\_$}
%\put(-4,0){$\_$} \put(-1,0){$\_$}\put(-2,0){$\_$}\put(-3,0){$\_$}
%\put(-5,0){$\_$}\put(-6,0){$\_$}\put(-7,0){$\_$}\put(-8,0){$\_$}
%\put(-10,0){$\_$}\put(-9,0){$\_$}\put(-11,0){$\_$}\put(-12,0){$\_$}
\put(0,2){$\vdots$}\put(15,2){$\vdots$}
\put(-4,-12){$\ldots$}\put(-8,-12){$\ldots$}\put(-12,-12){$\ldots$}
\put(-4,0){$\ldots$} \put(-8,0){$\ldots$}\put(-12,0){$\ldots$}
\put(0,-3){\thinlines\line(1,0){15}}\put(0,-6){\thinlines\line(1,0){15}}
\put(0,-9){\thinlines\line(1,0){9}}\put(0,-12){\thinlines\line(1,0){6}}
\put(0, 0){\thinlines\line(0,-1){12}}
\put(-12,6){\vector(1,0){32}}\put(21, 6){X}
\put(-12,6){\vector(0,-1){32}}\put(-15, -28){Y}
 \end{picture}}
 \caption{A lattice path from $(a,b)$ to $(c,d)$ and its associated Ferrers diagram}
\end{figure}
%%%%%%%%%%%%%%%%%%%%%%%%%%%%%%%%%%%%%%%%%%%%%%%%%%%%%%%%%%%%

Two lattice paths are said to be \emph{disjoint} or \emph{non
intersection} if they have no common points. For each path $w$
from $A$ to $B$ with $l$ vertical steps of abscissa $x_1,
x_2,\ldots, x_l$, arranged in decreasing order,
 we can associate a partition of integers
$\lambda_w=(x_1-a, x_2-a, \ldots, x_l-a)$. Actually the Ferrers
graph of $\lambda_w$ corresponds to the area of the region limited
by the lines $x=a$, $y=d$ and the  horizontal and vertical steps
of $w$. The weight of the partition $\lambda_w$ is defined by
$$
|\lambda_w|=(x_1-a)+(x_2-a)+\cdots +(x_l-a).
$$
For example, for
the lattice path $w$ in Figure~2, we have
$|\lambda_w|=5+5+3+2=15$. Define the weight of a $n$-tuple
$\gamma=(\gamma_1,\gamma_2,\ldots, \gamma_n)$ of lattice paths by
$$
\psi(\gamma)=q^{|\lambda_{\gamma_1}|+\ldots +|\lambda_{\gamma_n}|
}.
$$

We need   the following result, which can be  easily verified.
\begin{lem}\label{detlem}
Let $(a_{ij})_{i,j=0,\ldots, m}$ be an invertible lower triangular
matrix, and let $(b_{ij})_{i,j}=(a_{ij})_{i,j}^{-1}$. Then for
$0\leq k\leq n\leq m$, we have $$ b_{n,k}=\frac{(-1)^{n-k}}
{a_{k,k}a_{k+1,k+1}\cdots a_{n,n}}\left|a_{k+i,
k+j-1}\right|_{i,j=1,\ldots, n-k}.
$$
\end{lem}

Let $\Gamma_{k,l}$ be the set of $n$-tuples of non intersecting
lattice paths $\gamma=(\gamma_1,\ldots, \gamma_n)$ such that
\begin{itemize}
\item $\gamma_i$ goes from $A_i(i-1,2i-1)$ to $B_i(2i-1,2i-1)$ for $1\leq i<l$ or $k<i\leq n$
and from $A_{i+1}(i,2i+1)$ to $B_i(2i-1,2i-1)$ for
$l\leq i<k$.
\end{itemize}
\begin{thm} For integers $k,\,l\geq 1$ the coefficient
$c_{k,\,l}(q)$ is the generating function of $\Gamma_{k,l}$ with
respect to the weight $\psi$, i.e.,
$$
c_{k,\,l}(q)=\sum_{\gamma \in \Gamma_{k,l}}q^{\psi(\gamma)}.
$$
\end{thm}
\pv By Lemma~\ref{detlem}, for $1\leq l\leq k$ and $n\geq k$, we
have
\begin{eqnarray*}\label{coeff}
c_{k,l}(q)&=&\left| {{l+i}\brack
{2i-2j+2}}_q q^{(i-j)(i-j+1)}\right|_{i,j=1}^{k-l}\\
&=&\left| {{l+i+1}\brack
{2i-2j+2}}_qq^{(i-j)(i-j+1)}\right|_{i,j=0}^{k-l-1}\\ &=&
\sum_{\sigma\in S_n}(-1)^{inv(\sigma)}\prod_{i=1}^n{{l+i+1}\brack
{2i-2\sigma(i)+2}}_q\,q^{(i-\sigma(i))(i-\sigma(i)+1)}.
\end{eqnarray*}
For any $\sigma\in S_n$ denote by $C(\sigma,k,l)$ the set of
$n$-tuples of lattice paths
 $\gamma=(\gamma_1,\cdots,\gamma_n)$, where $\gamma_i$ goes
from $A_i$ to $B_{\sigma(i)}$
  for $1\leq i<l$ or $k<i\leq n$,
and from $A_{i+1}$ to $B_{\sigma (i)}$ for $l\leq i<k$.

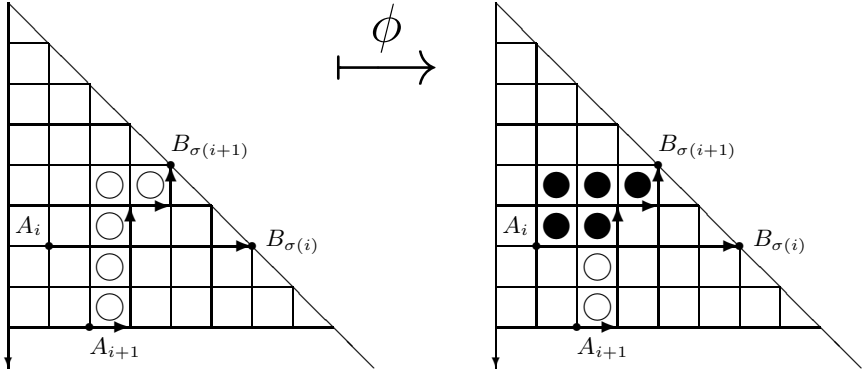
\begin{figure}
{\scriptsize
\setlength{\unitlength}{0.18cm}%
\begin{picture}(65,40)(-2,-37)
%\put(0,-3){\thicklines\thicklines\vector(1,0){3}}
\put(0,-18){\circle*{.5}}\put(3,-24){\circle*{.5}}\put(9,-12){\circle*{.5}}
\put(15,-18){\circle*{.5}} \put(0,-3){\thinlines\line(0,-1){21}}
\put(3, -6){\thinlines\line(0,-1){18}}
\put(6,-9){\thinlines\line(0,-1){15}}
\put(9,-12){\thinlines\line(0,-1){12}}
\put(-3,0){\thinlines\line(1,-1){27}}
\put(33,0){\thinlines\line(1,-1){27}} \put(12,
-15){\thinlines\line(0,-1){9}} \put(15,
-18){\thinlines\line(0,-1){6}} \put(18,
-21){\thinlines\line(0,-1){3}}
\put(-3,-3){\thinlines\line(1,0){3}}
\put(-3,-6){\thinlines\line(1,0){6}}
\put(-3,-9){\thinlines\line(1,0){9}}
\put(-3,-12){\thinlines\line(1,0){12}}
\put(-3,-15){\thinlines\line(1,0){15}}
\put(-3,-18){\thinlines\line(1,0){18}}
\put(-3,-21){\thinlines\line(1,0){21}}
\put(-3,-24){\thinlines\line(1,0){24}}
\put(-3,0){\vector(0,-1){27}} \put(-2.5,-17){$A_i$}
\put(3,-26){$A_{i+1}$} \put(9,-11){$B_{\sigma(i+1)}$}
\put(16,-18){$B_{\sigma(i)}$} \put(4.5,-13.5){\circle{2}}
\put(4.5,-16.5){\circle{2}}\put(4.5,-19.5){\circle{2}}
\put(4.5,-22.5){\circle{2}}\put(7.5,-13.5){\circle{2}}
\put(6,-15){\thicklines\thicklines\vector(1,0){3}}
\put(6,-24){\thicklines\thicklines\vector(0,1){9}}
\put(9,-15){\thicklines\thicklines\vector(0,1){3}}
\put(3,-24){\thicklines\thicklines\vector(1,0){3}}
\put(0,-18){\thicklines\thicklines\vector(1,0){15}}
\put(21,-6){\Huge {$\stackrel{\phi}{\longmapsto}$}}
\put(36,-18){\circle*{.5}}\put(39,-24){\circle*{.5}}\put(45,-12){\circle*{.5}}
\put(51,-18){\circle*{.5}} \put(36,-3){\thinlines\line(0,-1){21}}
\put(39, -6){\thinlines\line(0,-1){18}}
\put(42,-9){\thinlines\line(0,-1){15}}
\put(45,-12){\thinlines\line(0,-1){12}} \put(48,
-15){\thinlines\line(0,-1){9}} \put(51,
-18){\thinlines\line(0,-1){6}} \put(54,
-21){\thinlines\line(0,-1){3}}
\put(33,-3){\thinlines\line(1,0){3}}
\put(33,-6){\thinlines\line(1,0){6}}
\put(33,-9){\thinlines\line(1,0){9}}
\put(33,-12){\thinlines\line(1,0){12}}
\put(33,-15){\thinlines\line(1,0){15}}
\put(33,-18){\thinlines\line(1,0){18}}
\put(33,-21){\thinlines\line(1,0){21}}
\put(33,-24){\thinlines\line(1,0){24}}
\put(33,0){\vector(0,-1){27}} \put(33.5,-17){$A_i$}
\put(39,-26){$A_{i+1}$} \put(0,-30){}
\put(36,-30){}
\put(45,-11){$B_{\sigma(i+1)}$} \put(52,-18 ){$B_{\sigma(i)}$}
\put(40.5,-13.5){\circle*{2}}
\put(40.5,-16.5){\circle*{2}}\put(40.5,-19.5){\circle{2}}
\put(40.5,-22.5){\circle{2}}\put(43.5,-13.5){\circle*{2}}
\put(37.5,-13.5){\circle*{2}} \put(37.5,-16.5){\circle*{2}}
\put(42,-15){\thicklines\thicklines\vector(1,0){3}}
\put(42,-24){\thicklines\thicklines\vector(0,1){9}}
\put(45,-15){\thicklines\thicklines\vector(0,1){3}}
\put(39,-24){\thicklines\thicklines\vector(1,0){3}}
\put(36,-18){\thicklines\thicklines\vector(1,0){15}}
\end{picture}
\vspace{-1cm}
 \caption{\label{change}Change of weight after switching tails.}}
\end{figure}

Let $f:S_n\to \Z$ be a mapping defined by:
$$\forall \sigma\in S_n,
\quad f(\sigma)=\sum_{i=1}^n(i-\sigma(i))(i-\sigma(i)+1).
$$
Since the $q$-binomial coefficient has the following
interpretation~\cite[p. 33]{An}:
$$
{{m+n}\brack  m}_q=\sum_{\gamma}q^{|\lambda_\gamma|},
$$
where the sum is over all lattice  paths $\gamma$ from $(0,m)$ to
$(n,0)$, we derive immediately
\begin{equation}\label{eq:crucial}
c_{k,l}(q)=\sum_{\sigma\in S_n}\sum_{\gamma\in
C(\sigma,k,l)}(-1)^{inv(\sigma)}q^{\psi(\gamma)+f(\sigma)}.
\end{equation}
For any $n$-tuple of lattice paths $(\gamma_1,\ldots, \gamma_n)$,
if there is at least one intersecting point, we can define the
\emph{extreme intersecting point} $(i,j)\in\Pi$ to be  the
greatest intersecting point by the lexicographic order of their
coordinates. It is easy to see that this point must be an
intersecting point of two  lattice paths $w_i$ and $w_{i+1}$ of
consecutive indices. Applying the Gessel-Viennot method by
"switching the tails", i.e., exchanging the parts of $w_i$ and
$w_{i+1}$ starting from  the extreme point. Let $\phi:
\gamma\mapsto \gamma'$ be the corresponding  transformation on the
$n$-tuple of lattice paths with at least one intersecting point.
This transformation doesn't keep the value $\psi$ of intersecting
paths as illustrated in Figure~\ref{change}. However, it is easy
to see that $f$ is the unique mapping on $S_n$ satisfying
$f(id)=0$ and
 $$
 f(\sigma)-f(\sigma\circ
 (i,i+1))=2(\sigma(i)-\sigma(i+1)), \qquad \textrm{for any}\quad  \sigma\in S_n.
 $$
Hence, for any $\sigma\in S_n$ and $\gamma\in C({\sigma},k,l)$, we
have:
$$
q^{\psi(\gamma)+f(\sigma)}(-1)^{inv(\sigma)}=
 -q^{\psi(\phi(\gamma))+f(\sigma\circ (i,i+1))}
 (-1)^{\textrm{inv}(\sigma\circ (i,i+1))}.
$$
It means that $\phi$ is a \emph{weight-preserving-sign-reversing}
involution on the set of $n$-tuples of intersecting lattice paths
in $\cup_{\sigma \in S_n}C(\sigma,k,l)$. As $\gamma\in
C(\sigma,k,l)$ is non-intersecting only if $\sigma$ is an identity
permutation, that is $\gamma\in C(id,k,l)$. The result follows
then from \eqref{eq:crucial}.
 \qed
%%%%%%%%%%%%%%%%%%%%%%%%%%%%%%%%%%%%%%%%%%%%%%%%%%%%%%%%%%%%%%%%%%%%%%%%%

Notice that for $1\leq i<l$ or $k<i\leq n$, there is only one
lattice path from $A_i$ to $B_i$, the others have two vertical
steps. To each vertical step of $\gamma_i$ we can associate the
number $v=x_0-i+1$ between 1 and $i$, where $x_0$ is the abscissa
of the vertical step. We define the function $p:
[2n-2]\longrightarrow [0,\, n-1]$ as follows~:
$$
p(i)=\left\{\begin{array}{ll} 0 &
\textrm{if there is no vertical steps between the lines}\quad
y=i,\, y=i+1;\\
 v& \textrm{if v is the number associated to the
vertical step}
\end{array}\right.
$$
For example, for the preceding configuration, we have
$$
p(1)=\ldots =p(4)=0,\, p(5)=2,\, p(6)=1,\,p(7)=p(8)=p(10)=3,\,
p(9)=5.
$$
By construction, $p(2i-1)\geq p(2i)$ for all $i\in
[n-1]$. Now the condition of non-intersecting paths is equivalent
to $p(2i)\leq p(2i+1)$ for all $i\in [k-2]\setminus[l-1]$; and the
value of $w$ is $\psi(w)=-2(n-k)+\sum_ip(i)$.
%%%%%%%%%%%%%%%%%%%%%%%%%%%%%%%%%%%%%%%%%%%%%%%%%%%%%%%%%%%%%%%%%%%%%%%%%%%%%%%%%%%%%%%%%%%%%%%%%%%%%%%%%%%
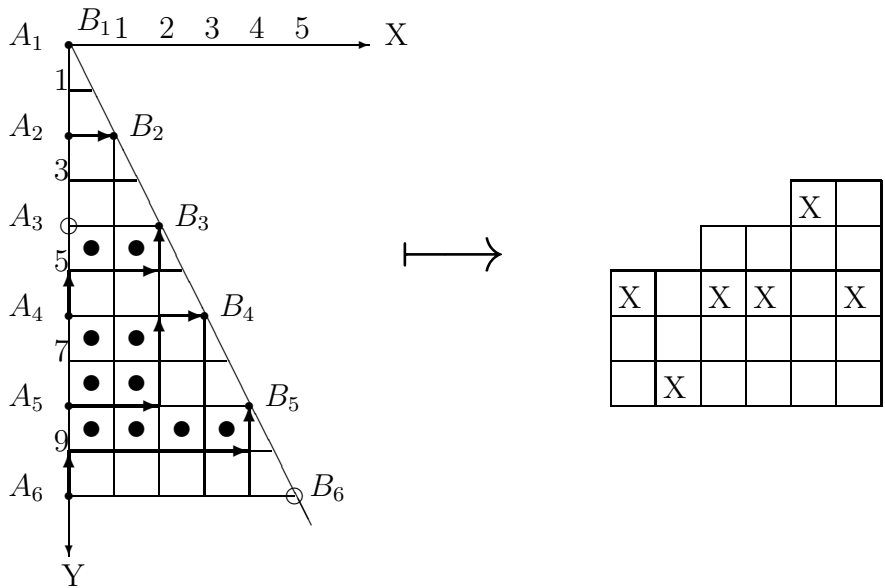
\begin{figure}
 \setlength{\unitlength}{0.20cm}
\begin{picture}(70,45)(-10,-40) \put(-4,0){$A_1$}\put(.5,1){$B_1$}\put(0,0){\circle*{.5}}
\put(0,0){\vector(1, 0){20}}\put(21, 0){X}
\put(3,.5){1}\put(6,.5){2}\put(9,.5){3}\put(12,.5){4}\put(15
,.5){5} \put(-1,-3){1}\put(-1,-6){}
\put(0,-6){\circle*{.5}}\put(-4,-6){$A_2$}
\put(0,-6){\thicklines\thicklines\vector(1,0){3}}
\put(3,-6){\circle*{.5}}\put(4,-6){$B_2$}
\put(0,-12){\circle{1}}\put(-4,-12){$A_3$}
\put(6,-12){\circle*{.5}}\put(7,-12){$B_3$}
\put(0,-18){\thicklines\thicklines\vector(0,1){3}}
\put(0,-15){\thicklines\thicklines\vector(1,0){6}}
\put(6,-15){\thicklines\thicklines\vector(0,1){3}}
\put(0,-18){\circle*{.5}}\put(-4,-18){$A_4$}
\put(9,-18){\circle*{.5}}\put(10,-18){$B_4$}
\put(0,-24){\circle*{.5}}\put(-4,-24){$A_5$}
\put(0,-24){\thicklines\thicklines\vector(1,0){6}}
\put(6,-24){\thicklines\thicklines\vector(0,1){6}}
\put(6,-18){\thicklines\thicklines\vector(1,0){3}}
\put(12,-24){\circle*{.5}} \put(1.5,-13.5){\circle*{1}}
\put(4.5,-13.5){\circle*{1}}\put(1.5,-19.5){\circle*{1}}
\put(4.5,-19.5){\circle*{1}}\put(1.5,-22.5){\circle*{1}}
\put(4.5,-22.5){\circle*{1}}
\put(1.5,-25.5){\circle*{1}}\put(4.5,-25.5){\circle*{1}}
\put(7.5,-25.5){\circle*{1}} \put(10.5,-25.5){\circle*{1}}
\put(13,-24){$B_5$} \put(0,-30){\circle*{0.5}}\put(-4,-30){$A_6$}
\put(0,-30){\thicklines\thicklines\vector(0,1){3}}
\put(0,-27){\thicklines\thicklines\vector(1,0){12}}
\put(12,-27){\thicklines\thicklines\vector(0,1){3}}
\put(15,-30){\circle{1}}\put(16,-30){$B_6$}
\put(-1,-9){3}\put(-1,-12){}\put(-1,-15){5}
\put(-1,-18){}\put(-1,-21){7}\put(-1,-24){}\put(-1,-27){9}\put(-2,-30){}
\put(3, -6){\thinlines\line(0,-1){24}} \put(6,
-12){\thinlines\line(0,-1){18}} \put(9,
-18){\thinlines\line(0,-1){12}} \put(12,
-24){\thinlines\line(0,-1){6}}
\put(.1,.1){\thinlines\line(1,-2){16}}
\put(0,-3){\thinlines\line(1,0){1.5}}\put(0,-6){\thinlines\line(1,0){3}}
\put(0,-9){\thinlines\line(1,0){4.5}}\put(0,-12){\thinlines\line(1,0){6}}
\put(0,-15){\thinlines\line(1,0){7.5}}
\put(0,-18){\thinlines\line(1,0){9}}\put(0,-21){\thinlines\line(1,0){10.5}}
\put(0,-24){\thinlines\line(1,0){12}}
\put(0,-27){\thinlines\line(1,0){13.5}}
\put(0,-30){\thinlines\line(1,0){15}}
\put(0,0){\vector(0,-1){34}}\put(-.5,-36){Y}\put(-3,-39){}
\put(22,-15){\Huge{$\longmapsto$}}
\put(36,-24){\thinlines\line(1,0){18}}
\put(36,-21){\thinlines\line(1,0){18}}
\put(36,-18){\thinlines\line(1,0){18}}
\put(36,-15){\thinlines\line(1,0){18}}
\put(42,-12){\thinlines\line(1,0){12}}
\put(48,-9){\thinlines\line(1,0){6}}
\put(54,-9){\thinlines\line(0,-1){15}}
\put(51,-9){\thinlines\line(0,-1){15}}
\put(48,-9){\thinlines\line(0,-1){15}}
\put(45,-12){\thinlines\line(0,-1){12}}
\put(42,-12){\thinlines\line(0,-1){12}}
\put(39,-15){\thinlines\line(0,-1){9}}
\put(36,-15){\thinlines\line(0,-1){9}}
\put(36.5,-17.5){X}\put(39.5,-23.5){X}\put(42.5,-17.5){X}\put(45.5,-17.5){X}\put(48.5,-11.5){X}\put(51.5,-17.5){X}
\end{picture}
\caption{\label{dpath}One of the 493 configurations counted by
$d_{6,3}(1)$ and its associated truncated pistol.}
\end{figure}
%%%%%%%%%%%%%%%%%%%%%%%%%%%%%%%%%%%%%%%%%%%%%%%%%%%%%%%%%%%%

Then we obtain a bijection between the configurations of
Proposition~5 and those that we can call \emph{truncated
alternating pistols}. More precisely we have the following result:
\begin{thm}
For $0\leq l\leq k$ and $n\geq k$, the coefficient
$c_{k+1,l+1}(q)$ is the generating function of alternating pistols
of $[2k]$, weighted by ${\ch'}$ and truncated at the index $2l$,
i.e. the weight of mappings $p: [2k]\longrightarrow [0,\, k]$
satisfying the three conditions:
\begin{enumerate}
\item $p(2i-1)=p(2i)=0$ for $1\leq i\leq l$,
\item $p(2i-1)\leq i$ and $p(2i)\leq i$ for $l<i\leq k$,
\item $p(2i-1) \geq p(2i)\leq p(2i+1)$ for $1\leq i<k$.
\end{enumerate}
\end{thm}

For example, the array $(g_{i,j}')$ with $5\leq i\leq 8$ and
$1\leq j\leq 4$, corresponding to the truncated alternating
pistols using for counting the coefficient
$c_{5,3}(q)=\sum_{k=1}^4 q^{k-1}g_{8,k}'$ is given in
Table~\ref{compuc53}.

%%%%%%%%%%%%%%%%%%%%%%%%%%%%%%%%%%%%%%%%%%%%%%%%%%%%%%%%%%%%%%%%%%%%%%%%%%%%%%%%%%%%%%%%%%%%%%%%%%%%%%%%%

\begin{table}[ht]{\scriptsize
\begin{tabular}{ccccccc|}
\multicolumn{7}{c}{ }\\ \multicolumn{7}{c}{ }\\
\multicolumn{7}{c}{ }\\ \multicolumn{7}{c}{ }\\
 \cline{5-7}
&&&&\multicolumn{1}{|c}{$1+q+2q^2+q^3+q^4$}&\multicolumn{1}{c|}
{$q^3+q^4+2q^5+q^6+q^7$}&4\\
 \cline{3-4}
 &&\multicolumn{1}{|c}{$1$}&\multicolumn{1}{c|}
 {$q^2$}&\multicolumn{1}{|c}{$1+q+2q^2+q^3+q^4$}&\multicolumn{1}{c|}
 {$q^2+2q^3+3q^4+3q^5+2q^6+q^7$}&3\\
&&\multicolumn{1}{|c}{$1$}&\multicolumn{1}{c|}{$q+q^2$}&\multicolumn{1}{|c}
{$1+q+2q^2+q^3$}&\multicolumn{1}{c|}
{$q+2q^2+4q^3+4q^4+3q^5+2q^6+q^7$}&2\\
&&\multicolumn{1}{|c}{$1$}&\multicolumn{1}{c|}
{$1+q+q^2$}&\multicolumn{1}{|c}{$1+q+q^2$}&\multicolumn{1}{c|}
{$1+2q+3q^2+4q^3+4q^4+3q^5+2q^6+q^7$}&1\\ \cline{3-7}
&&\multicolumn{1}{|c|}{$5$}&\multicolumn{1}{|c|}{$6$}
&\multicolumn{1}{|c|}{$7$}&\multicolumn{1}{|c|}{$8$}&$i\setminus
j$\\ \cline{3-7} \multicolumn{7}{c}{ }\\
\end{tabular}}
\caption{\label{compuc53}Computation of $c_{5,3}(q)$}
\end{table}
%%%%%%%%%%%%%%%%%%%%%%%%%%%%%%%%%%%%%%%%%%%%%%%%%%%%%%%%%%%%%%%%%%%%%%%%%%%%%%%%%%%
In particular we recover the alternating pistol in the case $l=0$,
and then we obtain the following result:

\begin{cor} For $n\geq 1$, the coefficient $c_{n,1}(q)$ of the inverse
matrix of $G_q$ is the $q$-Genocchi number $G_{2n}(q)$.
\end{cor}

 Now we give a last combinatorial interpretation of the
 $q$-Genocchi numbers. Some definitions about
 \emph{integer partitions} are needed.
 A partition $\mu=(\mu_1,\mu_2,\ldots)$ is said to \emph{smaller} than another partition
 $\lambda=(\lambda_1,\lambda_2, \ldots)$ if and only if all
the parts of $\mu$ are smaller than the one of $\lambda$. If
$\mu\leq \lambda$ we define a skew hook of shape $\lambda\setminus
\mu$ as the set difference of the diagram of $\lambda$ removed
that of $\mu$. Finally, a row-strict plane partition $T$ of
$\lambda\setminus \mu$ is a skew hook of shape $\lambda\setminus
\mu$ where we associate to the $j^{th}$ cell (from left to right)
of the $i^{th}$ line (from top to bottom), an positive integer
$p_{i,j}(T)$ such that, $\forall i\in [k],\,\forall j\in
[\lambda_i-\mu_i]$:
\begin{equation}\label{defpart}
p_{i,j}(T)>p_{i,j+1}(T)\quad\mbox{and}\quad p_{i,j}(T)\geq
p_{i+1,j}(T).
\end{equation}
 A reverse plane partition is obtained  by reversing all the
inequalities of (\ref{defpart}).

%%%%%%%%%%%%%%%%%%%%%%%%%%%%%%%%%%%%%%%%%%%%%%%%%%%%%%%%%%%%%%%%%%%%%%

\begin{figure}
\setlength{\unitlength}{0.15cm}
\begin{picture}(75,48)(-5,-45)
\put(-3.5,-3){$A_1$}\put(4,-3){$B_1$}
\put(0,-3){\thicklines\thicklines\vector(1,0){3}}
\put(0,-3){\circle*{.5}}
\put(10.5,-16.5){\circle*{1}}\put(13.5,-22.5){\circle*{1}}
\put(-.5,.5){0}\put(3,-3){\circle*{.5}} \put(0,0){\vector(1,
0){40}} \put(41,
0){X}\put(13.5,-25.5){\circle*{1}}\put(16.5,-22.5){\circle*{
1}}\put(16.5,-25.5){\circle*{1}}
\put(3,.5){1}\put(6,.5){2}\put(9,.5){3}\put(25.5,-28.5){\circle*{1}}\put(16.5,-28.5){\circle*{1}}
\put(12,.5){4}\put(15,.5){5}\put(16.5,-31.5){\circle*{1}}\put(19.5,-28.5){\circle*{1}}
\put(18,.5){6}\put(21,.5){7}\put(24,.5){8}\put(19.5,-31.5){\circle*{1}}\put(22.5,-28.5){\circle*{1}}
\put(27,.5){9}\put(30,.5){10}\put(-2,-6){2}
\put(3,-9){\circle*{.5}}\put(0,-11.5){$A_2$}
\put(3,-9){\thicklines\thicklines\vector(1,0){6}}
\put(9,-9){\circle*{.5}}\put(10,-9){$B_2$}
\put(6,-15){\circle{1}}\put(3.5,-17.5){$A_3$}
\put(15,-15){\circle*{.5}}\put(16,-15){$B_3$}
\put(9,-21){\circle*{.5}}\put(6.5,-23.5){$A_4$}
\put(9,-21){\thicklines\thicklines\vector(0,1){3}}
\put(9,-18){\thicklines\thicklines\vector(1,0){3}}
\put(12,-18){\thicklines\thicklines\vector(0,1){3}}
\put(12,-15){\thicklines\thicklines\vector(1,0){3}}
\put(21,-21){\circle*{.5}}\put(22,-21){$B_4$}
\put(12,-27){\circle*{.5}}\put(9.5,-29.5){$A_5$}
\put(12,-27){\thicklines\thicklines\vector(1,0){6}}
\put(18,-27){\thicklines\thicklines\vector(0,1){6}}
\put(18,-21){\thicklines\thicklines\vector(1,0){3}}
\put(12,-27){\circle*{.5}}\put(28,-27){$B_5$}
\put(27,-27){\circle*{.5}}\put(-3,-33){11}
\put(15,-33){\circle*{.5}}\put(12.5,-35.5){$A_6$}
\put(15,-33){\thicklines\thicklines\vector(1,0){6}}
\put(21,-33){\thicklines\thicklines\vector(0,1){3}}
\put(21,-30){\thicklines\thicklines\vector(1,0){6}}
\put(27,-30){\thicklines\thicklines\vector(0,1){3}}
\put(33,-33){\circle{1}}\put(34,-33){$ B_6$}
\put(-2,-9){3}\put(-2,-12){4}\put(-2,-15){5}
\put(-2,-18){6}\put(-2,-21){7}\put(-2,-24){8}\put(-2,-27){9}
\put(-3,-30){10} \put(6, -6){\thinlines\line(0,-1){27}} \put(3,
-3){\thinlines\line(0,-1){30}} \put(9,
-9){\thinlines\line(0,-1){24}} \put(12,
-12){\thinlines\line(0,-1){21}} \put(15,
-15){\thinlines\line(0,-1){18}} \put(18,
-18){\thinlines\line(0,-1){15}} \put(21,
-21){\thinlines\line(0,-1){12}} \put(24,
-24){\thinlines\line(0,-1){9}} \put(27,
-27){\thinlines\line(0,-1){6}} \put(30,
-30){\thinlines\line(0,-1){3}} \put(0,-3){\thinlines\line(1,0){3}}
\put(0,-6){\thinlines\line(1,0){6}}
\put(0,-9){\thinlines\line(1,0){9}}
\put(0,-3){\thinlines\line(1,-2){16}}
\put(0,-12){\thinlines\line(1,0){12}}
\put(0,-15){\thinlines\line(1,0){15}}
\put(0,-18){\thinlines\line(1,0){18}}
\put(0,-21){\thinlines\line(1,0){21}}
\put(0,-24){\thinlines\line(1,0){24}}
\put(0,-27){\thinlines\line(1,0){27}}
\put(0,-30){\thinlines\line(1,0){30}}
\put(0,-33){\thinlines\line(1,0){33}}
\put(0,0){\vector(0,-1){39}}\put(-.5,-41){Y}\put(5,-43){}
\put(0,0){\thinlines\line(1,-1){35}}
\put(40,-15){\Huge{$\longmapsto$}}
\put(54,-24){\thinlines\line(1,0){18}}
\put(54,-21){\thinlines\line(1,0){18}}
\put(54,-18){\thinlines\line(1,0){18}}
\put(54,-15){\thinlines\line(1,0){18}}
\put(60,-12){\thinlines\line(1,0){12}}
\put(66,-9){\thinlines\line(1,0){6}}
\put(72,-9){\thinlines\line(0,-1){15}}
\put(69,-9){\thinlines\line(0,-1){15}}
\put(66,-9){\thinlines\line(0,-1){15}}
\put(63,-12){\thinlines\line(0,-1){12}}
\put(60,-12){\thinlines\line(0,-1){12}}
\put(57,-15){\thinlines\line(0,-1){9}}
\put(54,-15){\thinlines\line(0,-1){9}}
\put(54.5,-20.5){X}\put(57.5,-23.5){X}
\put(60.5,-17.5){X}\put(63.5,-17.5){X}\put(66.5,-11.5){X}\put(69.5,-17.5){X}
\end{picture} \caption{\label{path1}One
of the 736 configurations counted by $c_{6,3}(1)$ and its
associated truncated pistol.}
\end{figure}
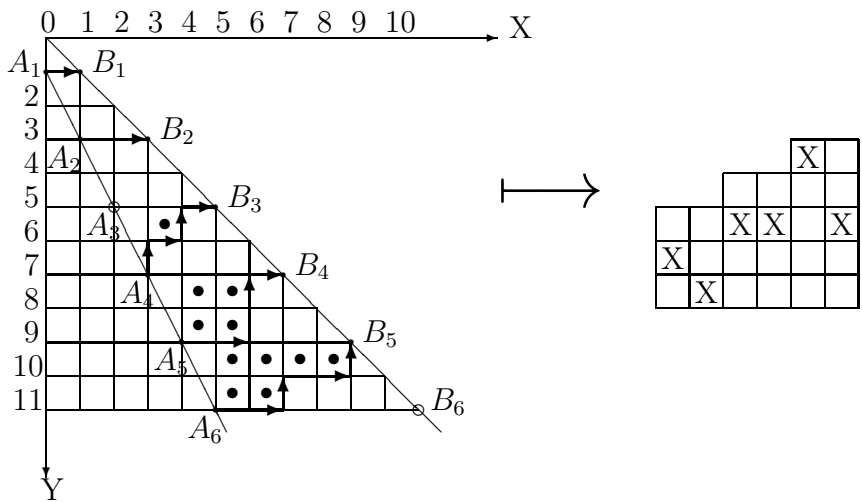
%%%%%%%%%%%%%%%%%%%%%%%%%%%%%%%%%%%%%%%%%%%%%%%%%%%%%%%%%%%%%%%%%%%%%%%%
Now, let $\gamma=(\gamma_1,\ldots, \gamma_n)$ be one of the
configuration counted by $c_{k,l}(q)$, $n\geq k\geq l$. Then we
can associate to this configuration, two partitions
$\lambda=(\lambda_1,\cdots,\lambda_n)$ and
$\mu=(\mu_1,\cdots,\mu_n)$ defined by $\lambda_i$ (resp. $\mu_i$)
equal $n+i-1$ for $i<l$ (resp. $i<k$) and $n+i+1$ otherwise. By
construction, $\lambda$ is larger  than $\mu$ and then we can
construct a row-strict plane partition $T$ where each case of
$\lambda\setminus \mu$ is labelled in the following way:

If the vertical steps of $\omega_{l+i-1}$ ($1\leq i\leq k-l$) have
 $x_{i,1}$ and $x_{i,2}$ for abscissa from left-to-right, so $x_{i,1}\leq
x_{i,2}$, define
$$
p_{i,j}(T)=2l+2i-j-x_{i,j}\qquad \hbox{for}\quad  j=1,2.
$$

For example, the row-strict plane partition corresponding to the
configuration of 5 paths in Figure~\ref{path1} is
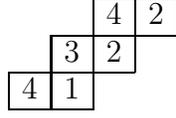
$$
\begin{tabular}{cccc}
\cline{3-4} & &\multicolumn{1}{|c|}4&\multicolumn{1}{|c|}2\\
\cline{2-4} &\multicolumn{1}{|c|}3&\multicolumn{1}{|c|} 2&\\
\cline{1-3} \multicolumn{1}{|c|}4&\multicolumn{1}{|c|}1&&\\
\cline{1-2}
\end{tabular}
$$

Let $T_{k,l}$ be the set of row-strict plane partition of form
$(k-l+1, k-l,\ldots, 2)-(k-l-1, k-l-2,\ldots, 0)$ such that the
largest entry in row $i$ is at most $l+i$. For any $T\in T_{k,l}$
define the value of $T$  by:
$$
|T|=\sum_{i=1}^{k-l}(p_{i,1}(T)+p_{i,2}(T)),
$$
then we have the following result, which is a $q$-analog of a
result of Gessel-Viennot~\cite[Theorem 31]{GV2}.

\begin{thm}
For $k\geq l\geq 1$, the entry $c_{k,l}(q)$ is the following
generating function of of $T_{k,l}$:
$$
c_{k,l}(q)=\sum_{T\in T_{k,l}}q^{k^2-l^2-|T|}.
$$
\end{thm}

%%%%%%%%%%%%%%%%%%%%%%%%%%%%%%%%%%%%%%%%%%%%%%%%%%%%%%%%%%%%%%%%%%%%%%%%
\section{Extension to negative indices and median $q$-Genocchi numbers}
%%%%%%%%%%%%%%%%%%%%%%%%%%%%%%%%%%%%%%%%%%%%%%%%%%%%%%%%%%%%%%%%%%%%%%%%
As in \cite{DZ2}, we can extend the matrix $G_q$ to the negative
indices  as follows~:
$$ H_q=\left({-j\brack 2i-2j}_qq^{(i-j)(2i-1)
}\right)_{i,j\geq 1}= \left({2i-j-1\brack j-1}_q \right)_{i,j\geq
1}, $$ and its inverse $$
H_q^{-1}=\left((-1)^{i-j}d_{i,j}(q)\right)_{i,j\geq 1}. $$

 Using the
result of Lemma~2, for $1\leq l\leq k$ and $n\geq k$, the
coefficient $d_{k,l}(q)$ is equal to:
\begin{equation}\label{det}
d_{k,l}(q)=\left|
{{l+2i-j}\brack {2i-2j+2}}_q\right|_{i,j=1}^{k-l}.
\end{equation}
The first values of $d_{i,j}(q)$ are given in Table~\ref{firdij}.

%%%%%%%%%%%%%%%%%%%%%%%%%%%%%%%%%%%%%%%%%%%%%%%%%%%%%%%%%%%%%
\begin{table}[ht] {\scriptsize
\begin {tabular}{cccc||c|c|c|c}
\multicolumn{8}{c}{ }\\ \multicolumn{8}{c}{ }\\
\multicolumn{8}{c}{ }\\ \multicolumn{8}{c}{ }\\&&&$i\setminus
j$&1&2&3&4\\ \cline{4-8} &&&1&1&0&0&0\\ &&&2&1&1&0&0\\
&&&3&${q}^{2}+q$ &${q}^{2}+q+1$&1&0\\
&&&4&$q^6+2q^5+2q^4+2q^3+q^2$&$q^6+2q^5+3q^4+3q^3+3q^2+q$&$
\left({q}^{2}+q+1\right )\left ({q}^{2}+1\right)$&1\\
\multicolumn{8}{c}{ }\\ \multicolumn{8}{c}{ }\\
\end{tabular}}
\caption{\label{firdij} First values of $d_{i,j}(q)$ }
\end{table}
%%%%%%%%%%%%%%%%%%%%%%%%%%%%%%%%%%%%%%%%%%%%%%%%%%%%%%%%%%%%%%%
 As in the previous section,
we then derive from \eqref{det} the following result.
\begin{thm} For integers $k,l\geq 1$ the coefficient
$d_{k,l}(q)$ is the generating function of configuration of
lattice path $\Omega=(\omega_1,\ldots, \omega_n)$, weighted by
$\psi$, satisfying the following two conditions :
\begin{enumerate}
\item $\omega_i$ joins $A_i(0,2i-2)$ to $B_i(i-1,2i-2)$ for $1\leq i<l$ or $k<i\leq n$
and $\omega_i$ joins $A_{i+1}(0,2i)$ to $B_i(i-1,2i-2)$ for $l\leq
i<k$.
\item the paths $\omega_1,\ldots, \omega_n$ are disjoint.
\end{enumerate}
\end{thm}

Similarly to the preceding section, remark that for $1\leq i<l$ or
$k<i\leq n$, there is an only lattice path from $A_i$ to $B_i$ and
the other ones have two vertical steps. To each vertical steps of
$\omega_i$, we associate a number $v=x_0+1$ between 1 and $i$
where $x_0$ is the abscissa of this vertical step. Then we can
define a function $p: [2n-2]\longrightarrow [0,\, n-1]$  as
follows~:
 $$ p(i)=\left\{\begin{array}{ll} 0 & \textrm{if there is
no vertical steps between the lines}\quad y=i-1,\, y=i,\\ v&
\textrm{if v is the number associated to the vertical step}.
\end{array}
\right.
$$

For example, for the preceding configuration, we have
$p(1)=p(2)=p(3)=p(4)=0$, $p(5)=p(7)=p(8)=3$, $p(6)=p(10)=1$,
$p(9)=5$. By construction, $p(2i-1)\geq p(2i)$ for all $i\in
[n-1]$ and the condition of non-intersecting paths is equivalent
to $p(2i)< p(2i+1)$ for all $i\in [k-2]\setminus [l-1]$. The value
of $w$ is $\psi(w)=-2(n-k)+\sum_ip(i).$ Then we obtain a bijection
between the configurations of Proposition 8 and those that we can
call \emph{truncated alternating pistols}. More precisely we state
the following result:
\begin{prop}
For $0\leq l\leq k$ and $n\geq k$, the coefficient
$d_{k+1,l+1}(q)$ is the generating function of alternating pistols
of $[2k]$, weighted by $\ch'$ and truncated at the index $2l$,
i.e. the mappings $p: [2k]\longrightarrow [0,\, k]$ satisfying the
three conditions :
\begin{enumerate}
\item $p(2i-1)=p(2i)=0$ for $1\leq i\leq l$,
\item $p(2i-1)\leq i$ and $p(2i)\leq i$ for $l<i\leq k$,
\item $p(2i-1) \geq p(2i)< p(2i+1)$ for $1\leq i<k$.
\end{enumerate}
\end{prop}

The array  for the computation of $d_{5,3}(q)$ is given in
Table~\ref{compud53}.
%%%%%%%%%%%%%%%%%%%%%%%%%%%%%%%%%%%%%%%%%

\begin{table}[ht]{\scriptsize
\begin{tabular}{ccccccc|}
\multicolumn{7}{c}{ }\\ \multicolumn{7}{c}{ }\\
\multicolumn{7}{c}{ }\\ \multicolumn{7}{c}{ }\\
 \cline{5-7}
&&&&\multicolumn{1}{|c}{$1+q+2q^2+q^3+q^4$}&\multicolumn{1}{c|}
{$q^3+q^4+2q^5+q^6+q^7$}&4\\
 \cline{3-4}
 &&\multicolumn{1}{|c}{$1$}&\multicolumn{1}{c|}
 {$q^2$}&\multicolumn{1}{|c}{$1+q+2q^2+q^3$}&\multicolumn{1}{c|}
 {$q^2+2q^3+3q^4+3q^5+q^6+q^7$}&3\\
&&\multicolumn{1}{|c}{$1$}&\multicolumn{1}{c|}{$q+q^2$}&\multicolumn{1}{|c}
{$1+q+q^2$}&\multicolumn{1}{c|}
{$q+2q^2+3q^3+3q^4+3q^5+q^6+q^7$}&2\\
&&\multicolumn{1}{|c}{$1$}&\multicolumn{1}{c|}
{$1+q+q^2$}&\multicolumn{1}{|c}{$0$}&\multicolumn{1}{c|}
{$q+2q^2+3q^3+3q^4+3q^5+q^6+q^7$}&1\\ \cline{3-7}
&&\multicolumn{1}{|c|}{$5$}&\multicolumn{1}{|c|}{$6$}
&\multicolumn{1}{|c|}{$7$}&\multicolumn{1}{|c|}{$8$}&$i\setminus
j$\\ \cline{3-7} \multicolumn{7}{c}{ }\\
\end{tabular}}
\caption{\label{compud53}Computation of $d_{5,3}(q)$}
\end{table}
%%%%%%%%%%%%%%%%%%%%%%%%%%%%%%%%%%%%%%%%%%%%%%%%%%%%%%%%

In particular we recover the alternating pistol when $l=0$, and
then we obtain the following result:

\begin{cor} For $n\geq 1$, the coefficient $d_{n,1}(q)$ of the inverse
matrix of $H_q$ is the medians $q$-Genocchi number $H_{2n+1}(q)$.
\end{cor}

 Now, let $\Omega=(\omega_1,\ldots, \omega_n)$ be one of
the configuration counting by $d_{k,l}(1)$, $n\geq k\geq l$. Then
we can associate to this configuration, two partitions
$\lambda=(\lambda_1,\cdots,\lambda_n)$ and
$\mu=(\mu_1,\cdots,\mu_n)$ defined by $\lambda_i$ (resp. $\mu_i$)
equal $n+i-2$ for $i<l$ (resp. $i<k$) and $n+i$ otherwise. By
construction, $\lambda$ is bigger than $\mu$ and then we can
construct an array $T$ where each case of $\lambda\setminus \mu$
is labelled in the following way:

If the vertical steps of $\omega_{l+i-1}$ ($1\leq i\leq k-l$) have
respectively $x_{i,1}$ and $x_{i,2}$ for abscissa, ($x_{i,1}\leq
x_{i,2}$), then $p_{i,j}(T)=x_{i,j}+1$  for $j=1,2$.

For example the row-strict plane partition
corresponding to the configuration of 5 paths in Figure~\ref{dpath} is

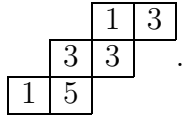
$$
\begin{tabular}{cccc}
\cline{3-4} & &\multicolumn{1}{|c|}1&\multicolumn{1}{|c|}3\\
\cline{2-4} &\multicolumn{1}{|c|}3&\multicolumn{1}{|c|} 3&\\
\cline{1-3} \multicolumn{1}{|c|}1&\multicolumn{1}{|c|}5&&\\
\cline{1-2}
\end{tabular}.
$$
Similarly  we have the following
\begin{thm}
For $k\geq l\geq 1$, $$d_{k,l}(q)=\sum_{T\in
\widetilde{T}_{k,l}}q^{-2(k-l)+|T|},$$ where $\widetilde{T}_{k,l}$
is the set of column-strict reverse plane partition of $(k-l+1,
k-l,\ldots, 2)-(k-l-1, k-l-2,\ldots, 0)$ with positive integer
entries in which the largest entry in row i is at most $l+i-1$.
\end{thm}

%%%%%%%%%%%%%%%%%%%%%%%%%%%%%%%%%%%%%%%%%%%%%%%%%%%%%%%%%%%%%%%%%%%%%%%%%%%%
\section{A remarkable triangle of $q$-numbers refining $q$-Euler numbers}
%%%%%%%%%%%%%%%%%%%%%%%%%%%%%%%%%%%%%%%%%%%%%%%%%%%%%%%%%%%%%%%%%%%%%%%%%%%%
Recall that the Euler numbers $E_{2n}$ are the coefficients in the
Taylor expansion of the function $\frac{1}{\cos x}$:
$$
\frac{1}{\cos x}=\sum_{n\geq 0}E_{2n}\frac{x^{2n}}{(2n)!}.
$$
Let $c_{i,j}=c_{i,j}(1)$. Then Dumont and Zeng~\cite{DZ1} proved
that there is a triangle of positive integers $k_{n,j}$ ($1\leq
j\leq n-1$) featuring the two kinds of Genocchi numbers and
refining Euler numbers as follows:
$$
k_{n,1}+k_{n,2}+\ldots k_{n,n-1}=E_{2n-2}, \quad k_{n,1}=G_{2n}
\quad\hbox{and}\quad  k_{n,n-1}=H_{2n-1}.
$$
Moreover,
$$
\sum_{j\geq 0}c_{n+j,j+1}x^{j+1}=\frac{k_{n,1}x+k_{n,2}x^2+\ldots
+k_{n,n-1}x^{n-1}}{(1-x)^{2n-1}}.
$$
The first values of $k_{n,j}$ ($1\leq j\leq n-1$) are tabulated as
follows:
$$
\begin{array}{c|cccccc|c}
n\setminus j&1&2&3&4&5&&\sum_{j}k_{n,j}=E_{2n-2}\\
\hline
1&1&&&&&&1\\
2&1&&&&&&1\\
3&3&2&&&&&5\\
4&17&36&8&&&&61\\
5&155&678&496&56&&&1385\\
6&2073&15820&23576&8444&608&&50521\\
\end{array}
$$
 We show now there is a $q$-analog of the above triangle.
Following Jackson~\cite{Ja} the $q$-secant numbers $E_{2n}(q)$ are
 defined by
$$\sum_{n\geq 0}E_{2n}(q)\frac{u^{2n}}{(q;q)_{2n}}=
\left(\sum_{n\geq 0}(-1)^n\frac{u^{2n}}{(q;q)_{2n}}\right)^{-1}.
$$

Let $[x]={(q^x-1)/(q-1)}$ and $[x]_n=[x][x-1]\cdots [x-n+1]$ for
$n\geq 0$. Then $([x]_n)$ is a basis of $C[q^x]$. For any integer
$n\geq 0$ we define a linear $q$-difference operator $\delta_q^n$
on $C[q^x]$ as follows~: for $f(x)\in C[q^x]$,
\begin{equation}
\delta_q^0f(x)=f(x),\qquad \delta_q^{n+1}f(x)=(E-q^nI)\,\delta_q^n
f(x).
\end{equation}
that is,
$$
\delta_q^nf(x)=(E-q^{n-1}I)(E-q^{n-2}I)\cdots (E-I)f(x).
$$
In view of the $q$-binomial formula~\cite[p. 36]{An}:
\begin{equation}\label{eq:qbin}
(x;q)_n=\sum_{k=0}^n(-1)^k{n\brack k}_qq^{k\choose 2}x^k,
\end{equation}
we have
$$
\delta_q^nf(x)=\sum_{k=0}^n(-1)^k{n\brack k}_q q^{k\choose
2}f(x+n-k).
$$

\begin{lem} For all non negative integers $n,m$ we have
$$ \delta_q^n[x]_m=\left\lbrace \begin{array}{ll}
[m]_n[x]_{m-n}q^{n(x+n-m)}&\quad \hbox{if}\quad n\leq m\\ 0&\quad
\hbox{if}\quad n>m.
\end{array}
\right. $$
\end{lem}
Hence $\delta_q^nf(x)=0$ if $f(x)$ is a polynomial in $q^x$ of
degree $<n$. It follows from the $q$-binomial
identity~\ref{eq:qbin} that
\begin{eqnarray*}
(x; q)_{2n-1} \sum_{j\geq 0}c_{n+j,j+1}(q)x^{j+1}
&=&\sum_{m\geq 0}x^{m+1}\sum_{k\geq 0}(-1)^k{{2n-1}\brack
k}_q q^{{k\choose 2}}c_{n+m-k,m-k+1}(q),\\ &=& \sum_{m\geq
0}x^{m+1}\delta_q^{2n-1}f(m).
\end{eqnarray*}
where $f(m)$ denotes the following determinant~:
$$f(m)=\left|{m-2(n-1)+i\brack
2i-2j+2}_qq^{(i-j)(i-j+1)}\right|_{i,j=1}^{n-1} $$ is a polynomial
in $q^m$ of degree $2(n-1)$ when $m\geq 2n-3$.
Hence the preceding
expression is a polynomial in $x$ of degree $d\leq 2n-1$, i.e., we have
 \begin{equation}\label{eq1}
\sum_{j\geq 0}c_{n+j,j+1}(q)x^{j+1}=\frac{\alpha_0(q)+\cdots+\alpha_{d-1}(q)
x^{d}}{(x; q)_{2n-1}}
 \end{equation}
Applying a well-known result about rational functions~\cite[p.
202-210]{St1}, we derive from (\ref{eq1}) that
\begin{eqnarray*}
\sum_{j\geq
1}c_{n-j,-j+1}(q)x^j&=&-\frac{\alpha_0+
\alpha_1 x^{-1}+\cdots+\alpha_{d-1}x^{-d}}
{(1/x;q)_{2n-2}}\\
&=&-\frac{\alpha_0x^{2n-1}+\cdots+\alpha_{d-1}x^{2n-d}}
{(x;q)_{2n-2}}.
\end{eqnarray*}
But the coefficient $c_{n-j,-j+1}(q)$ is null for all $1\leq j\leq
n$ because the determinant formula of  $c_{k,l}(q)$ contains a row
with only zeros.  So $d\leq n-1$.

Summarizing all the above we get the following theorem, which is a
$q$-analog of a result of Dumont and Zeng~\cite[Prop. 7]{DZ2}.
\begin{thm}\label{qtri}
 For $n\geq 2$, $\forall j\in [n-1]$, there
are polynomials $k_{n,j}(q)$ in $q$ such that \begin{eqnarray}
\sum_{j\geq 0}c_{n+j,j+1}(q)x^{j+1}&=&
\frac{\sum_{i=1}^{n-1}q^{(i-1)i}k_{n,i}(q)x^i}{(x;
q)_{2n-1}}.\label{kdef3}\\ \sum_{j\geq
0}d_{n+j,j+1}(q)x^{j+1}&=&\frac{\sum_{i=1}^{n-1}
q^{(i-1)i}k_{n,n-i}(q)x^{i}}{(x; q)_{2n-1}}.\label{kdef4}
\end{eqnarray}
Moreover, we
have
$k_{n,1}(q)=G_{2n}(q),\,\,k_{n,n-1}(q)=H_{2n-1}(q)$ and
$$
E_{2n-2}(q)=\sum_{i=1}^{n-1}{q^{(i-1)i}k_{n,n-i}(q)}.
$$
\end{thm}

\pv Equations (\ref{kdef3}) and (\ref{kdef4}) have been proved
previously. In view  of Corollaries~1 and 2 we derive from
(\ref{kdef3}) and (\ref{kdef4}) that
\begin{eqnarray*}
k_{n,1}(q)&=& c_{n,1}(q)=G_{2n}(q),\\
k_{n,n-1}(q)&=& d_{n,1}(q)=H_{2n-1}(q).
\end{eqnarray*}
Recall that for any sequence $(a_n)_n$ in $\C[[q]]$, we have
$\lim_{q\to 1}(1-x)\sum_{n\geq 0}a_nq^n=\lim_{n\to \infty}a_n$,
provided the later limit exists. Hence we derive from
(\ref{kdef4}) that
\begin{eqnarray*}
\sum_{i=1}^{n-1}{q^{(i-1)i}k_{n,n-i}(q)}&= &\lim_{x\rightarrow
1}(x; q)_{2n-1}\sum_{j\geq 0}d_{n+j,j+1}(q)x^{j+1}\\
&=& (q;q)_{2n-2}\;\lim_{j\to \infty}d_{n+j,j+1}(q).
\end{eqnarray*}
As $\lim_{n\to +\infty}{n\brack k}_q=\frac{1}{(q;q)_k}$ it follows
from (\ref{det}) that
\begin{equation}\label{eq:stan}
\sum_{i=1}^{n-1}{q^{(i-1)i}k_{n,n-i}(q)}=(q,q)_{2n-2}\left|
\frac{1}{(q;q)_{2i-2j+2}}\right|_{i,j=1}^{n-1}.
\end{equation}
Now, using  inclusion-exclusion principle we can show
(see~\cite[p.70]{St1}) that  the right-hand side of
\eqref{eq:stan} is the enumerating polynomial of up-down
permutations on $[2n-2]$, i.e., whose descent set is
$\{2,4,\cdots,2n-4\}$,  with respect to inversion numbers, and it
is also known (see \cite[p.148]{St1}) that  this enumerating
polynomial is equal to the $q$-Euler polynomial $E_{2n-2k}(q)$.
\qed

It is not difficult to derive from Theorem~\ref{qtri} the
following result.
\begin{cor}
 For $n\geq 2$, for all $i\in [n-1]$, we have:
$$ q^{(i-1)i}k_{n,i}(q)=\sum_{l=0}^{i-1}(-1)^lq^{{k\choose
2}}{2n-1\brack l}_qc_{n+i-l-1,\, i-l}(q), $$ and $$
q^{(i-1)i}k_{n,n-i}(q)=\sum_{l=0}^{i-1}(-1)^lq^{{k\choose
2}}{2n-1\brack l}_qd_{n+i-l-1,\, i-l}(q). $$
\end{cor}

Finally, for $n=2, \, 3$, equation (\ref{kdef3}) reads as follows:
\begin{eqnarray*}
\frac{x}{(x;q)_3}&=&x+(1+q+q^2)x^2+(1+q+2q^2+q^3+q^4)x^3+\cdots,\\
\frac{(1+q+q^2)x+q^2(q+q^2)x^2}{(x;q)_5}&=&(1+q+q^2)x\\
&+&(1+2q+3q^2+4q^3+4q^4+2q^5+q^6)x^2+\cdots.
\end{eqnarray*}
So $k_{3,1}(q)=1+q+q^2$ and $k_{3,2}(q)=q+q^2$. While the five
up-down permutations on $[4]$ are
$$
1\,3\,2\,4,\quad 1\,4\,2\,3,\quad 2\,3\,1\,4,\quad
2\,3\,1\,4,\quad 3\,4\,1\,2.
$$
Therefore $E_4(q)=q+2q^2+q^3+q^4$ and we can check that
$E_4(q)=k_{3,2}(q)+q^2k_{3,1}(q)$.

 For $n=4$ the
values of $k_{4,j}(q)$, $1\leq j\leq 3$, are given by
\begin{eqnarray*}
k_{4,1}(q)&=& 1+2q+3q^2+4q^3+4q^4+2q^5+q^6,\\
k_{4,2}(q)&=&q(1+q)(1+q^2)(1+q+q^2)^2,\\
k_{4,3}(q)&=&q^2(q^2+1)(q+1)^2.
\end{eqnarray*}

It seems that the coefficients of the polynomial $k_{n,i}(q)$ in
$q$ are \emph{non negative integers} and it would be interesting
to find a combinatorial interpretation for $k_{n,i}(q)$ in case
the above conjecture is true.

\end{document}